\documentclass[10pt]{article}
\usepackage{amsfonts}
\usepackage{hyperref}
\usepackage{mathrsfs,enumerate,amssymb,amsmath,color,lineno}
\usepackage{indentfirst}
\usepackage{amsmath}
\usepackage{amssymb}
\usepackage[amsmath,thmmarks]{ntheorem}
\usepackage{cite}
\usepackage{graphicx}
\usepackage{tikz}

\newtheorem{definition}{\bf Definition}[section]
\newtheorem{lemma}{\bf Lemma}[section]
\newtheorem{theorem}{\bf Theorem}[section]
\newtheorem{remark}{\bf Remark}[section]
\newtheorem{corollary}{\bf Corollary}[section]
\newtheorem{example}{\bf Example}[section]

\newtheorem{proposition}{\bf Proposition}[section]

\begin{document}
\setcounter{page}{1}

\title{{\textbf{The characterization of monotone functions that generate associative functions}}\thanks {Supported by
the National Natural Science Foundation of China (Nos.12071325, 12471440)}}
\author{Meng Chen\footnote{\emph{E-mail address}: mathchen2019@163.com}, Yun-Mao Zhang\footnote{\emph{E-mail address}: 2115634219@qq.com}, Xue-ping Wang\footnote{Corresponding author. xpwang1@hotmail.com; fax: +86-28-84761502},\\
\emph{School of Mathematical Sciences, Sichuan Normal University,}\\
\emph{Chengdu 610066, Sichuan, People's Republic of China}}

\newcommand{\pp}[2]{\frac{\partial #1}{\partial #2}}
\date{}
\maketitle
\begin{quote}
{\bf Abstract} Consider a two-place function $T: [0,1]^2\rightarrow [0,1]$ defined by $T(x,y)=f^{(-1)}(F(f(x),f(y)))$ where $F:[0,\infty]^2\rightarrow[0,\infty]$ is an associative function, $f: [0,1]\rightarrow [0,\infty]$ is a monotone function that satisfies either $f(x)=f(x^{+})$ when $f(x^{+})\in \mbox{Ran}(f)$ or $f(x)\neq f(y)$ for any $y\neq x$ when $f(x^{+})\notin \mbox{Ran}(f)$ for all $x\in[0,1]$ and $f^{(-1)}:[0,\infty]\rightarrow[0,1]$ is a pseudo-inverse of $f$. In this article, the associativity of the function $T$ is shown to depend only on properties of the range of $f$. The necessary and sufficient conditions for the $T$ being associative are presented by applying the properties of the monotone  function $f$.

\end{quote}

\section{Introduction}
The construction methods of associative functions starting from an one-place function play an indispensable role in the theory of solving associative equations \cite{CA2006}. The idea that a strictly monotone function may generate an associative function traces to the old work of Abel \cite{Abel}. The relation between triangular norms (t-norms for short) and continuous strictly decreasing additive generators was originally explored by Schweizer and Sklar \cite{BS1961} and Ling \cite{CH1965}, respectively. Note that the concept of an additive generator, which was first introduced by Schweizer and Sklar \cite{BS1961}, of a t-norm has been generalized several times in the literature, see e.g., \cite{SJ1999, EP2000,CH1965, AM2004, BS1963, YO2007, YO2008, PV2005, PV2008, DZ2005}. In particular, Klement, Mesiar and Pap \cite{EP2000} defined an additive generator as a strictly decreasing function $f :[0,1]\rightarrow [0,\infty]$ that is right-continuous at $0$ with $f (1) = 0$ such that for all $(x,y)\in [0,1]^2$,
\begin{equation}
\label{eq:1}
f(x)+f(y)\in \mbox{Ran}(f)\cup [f(0),\infty],
\end{equation}
and
\begin{equation}
\label{eq:2}
 T(x,y)= f^{(-1)}(f(x)+f(y))
\end{equation}
where $ f^{(-1)}$ is a pseudo-inverse of $f$ and $\mbox{Ran}(f)$ is a range of $f$. Based on Vicen\'{\i}k \cite{PV1998,PV1998b} they claimed that there exists a generalized additive generator just satisfying \eqref{eq:2}, and this was identified by Vicen\'{\i}k \cite{PV2005} when $f$ is a strictly monotone function and by Zhang and Wang \cite{YM2024} when $f$ is a monotone right continuous function, respectively. Yao Ouyang et al. \cite{YO2008} also generalized the concept of an additive generator of a t-norm as that there are a strictly decreasing function $f: [0,1]\rightarrow [0,\infty]$ with $f (1) = 0$ and a binary operator $\star: [0,\infty]^2\rightarrow [0,\infty]$ satisfying that $([0,\infty],\star,\leq)$ is a fully ordered Abel semigroup with neutral element $0$, such that
\begin{equation}
\label{eq:3}
 f(x)\star f(y)\in \mbox{Ran}(f)\cup [f(0^+),\infty]
\end{equation}
for all $x,y \in [0,1]$ and
\begin{equation}
\label{eq:4}
 T(x,y)= f^{(-1)}(f(x)\star f(y))
\end{equation}
where $ f^{(-1)}$ is the pseudo-inverse of $f$ and $\mbox{Ran}(f)$ is the range of $f$. Motivated by the work of Vicen\'{\i}k \cite{PV2005}, Yao Ouyang et al. \cite{YO2008} and Zhang and Wang \cite{YM2024}, this article further considers the problem: what is a characterization of monotone functions $f: [0,1]\rightarrow [0,\infty]$ such that the function $T: [0,1]^2\rightarrow [0,1]$ given by
\begin{equation}\label{eq:5}
T(x,y)= f^{(-1)}(F(f(x),f(y)))
\end{equation}
is associative where $f^{(-1)}$ is the pseudo-inverse of $f$ and $F: [0,\infty]^2\rightarrow [0,\infty]$ is an associative function?

The rest of this article is organized as follows. In Section 2, we show some preliminaries concerned on basic concepts and known results of t-norms and t-conorms, respectively. In Section 3, we give a representation of the range $\mbox{Ran}(f)$ of a kind of non-decreasing functions $f$. In Section 4, we first define an operation $\otimes$ on the $\mbox{Ran}(f)$, and then supply a necessary and sufficient condition for the operation $\otimes$ being associative. In Section 5, we characterize what properties of $\mbox{Ran}(f)$ are equivalent to the associativity of the function $T$ defined by Eq.(\ref{eq:5}). A conclusion is drawn in Section 6.

\section{Preliminaries}

This section recalls some known elementary concepts and results in the literature, respectively.
\begin{definition}[\cite{EP2000}]
\emph{A t-norm is a binary operator $T:[0, 1]^2\rightarrow [0, 1]$ such that for all $x, y, z\in[0, 1]$ the following conditions are satisfied:}

\emph{(T1)} $T(x,y)=T(y,x)$,

\emph{(T2)} $T(T(x,y),z)=T(x,T(y,z))$,

\emph{(T3)} $T(x,y)\leq T(x,z)$ \emph{whenever} $y\leq z$,

\emph{(T4)} $T(x,1)=x$.

\emph{A binary operator $T:[0, 1]^2\rightarrow [0, 1]$ is called a t-subnorm if it satisfies (T1), (T2), (T3) and $T(x,y)\leq \min\{x,y\}$ for all $x,y\in [0, 1]$.}
\end{definition}
\begin{definition}[\cite{EP2000}]
\emph{A t-conorm is a binary operator $S:[0, 1]^2\rightarrow [0, 1]$ such that for all $x, y, z\in[0, 1]$ the following conditions are satisfied:}

\emph{(S1)} $S(x,y)=S(y,x)$,

\emph{(S2)} $S(S(x,y),z)=S(x,S(y,z))$,

\emph{(S3)} $S(x,y)\leq S(x,z)$ \emph{whenever} $y\leq z$,

\emph{(S4)} $S(x,0)=x$.

\emph{A binary operator $S:[0, 1]^2\rightarrow [0, 1]$ is called a t-superconorm if it satisfies (S1), (S2), (S3) and $S(x,y)\geq \max\{x,y\}$ for all $x,y\in [0, 1]$.}
\end{definition}

\begin{definition}[\cite{EP2000,PV2005}]\label{def2.3}
\emph{Let $p, q, s, t\in [-\infty, \infty]$ with $p<q, s<t$ and $f:[p,q]\rightarrow[s,t]$ be a non-decreasing (resp. non-increasing) function. Then the function $f^{(-1)}:[s,t]\rightarrow[p,q]$ defined by
\begin{equation*}
f^{(-1)}(y)=\sup\{x\in [p,q]\mid f(x)<y\}\,(\mbox{resp. }f^{(-1)}(y)=\sup\{x\in [p,q]\mid f(x)>y\})
\end{equation*}
is called a pseudo-inverse of the non-decreasing (resp. non-increasing) function $f$.}
\end{definition}

Denoted by $A\setminus B=\{x\in A\mid x\notin B\}$ for two sets $A$ and $B$.

\section{The range of a non-decreasing  function}
This section describes the range of a kind of non-decreasing functions.

Let $f:[0,1]\rightarrow [0,\infty]$ be a function. We write $f(a^-)=\lim_{x \rightarrow a^{-}}f(x)$ for each $a\in (0, 1]$ and $f(a^+)=\lim_{x\rightarrow a^{+}}f(x)$ for each $a\in [0, 1)$. Define $f(0^-)=0$ and $f(1^+)=\infty$ whenever $f$ is non-decreasing.

Let $M\subseteq [0,\infty]$. A point $x\in [0,\infty]$ is said to be an accumulation point of $M$ from the left if there is a strictly increasing sequence $\{x_{n}\}_{n\in N}$ of points $x_{n}\in M$ satisfying $\lim_{n\rightarrow \infty}x_{n}=x$. A point $x\in [0,\infty]$ is said to be an accumulation point of $M$ from the right if there is a strictly decreasing sequence $\{x_{n}\}_{n\in N}$ of points $x_{n}\in M$ satisfying $\lim_{n\rightarrow \infty}x_{n}=x$. A point $x\in [0,1]$ is said to be an accumulation point of $M$ if it is an accumulation point of $M$ from the left or from the right. Denote the set of all accumulation points of $M$ by Acc$(M)$.

Denote by $\mathcal{F}$ the set of all non-decreasing functions $f:[0,1]\rightarrow [0,\infty]$ that satisfy either $f(x)=f(x^{+})$ when $f(x^{+})\in \mbox{Ran}(f)$ or $f(x)\neq f(y)$ for any $y\neq x$ when $f(x^{+})\notin \mbox{Ran}(f)$ for all $x\in[0,1]$. It is easy to see that both a strictly increasing and non-decreasing right continuous functions are special elements of $\mathcal{F}$.
Let
$$\mathcal{A}=\{M \mid \mbox{there is a function } f\in \mathcal{F}\mbox{ such that }\mbox{Ran}(f)=M\}.$$
Then the following lemma describes the range of a function $f\in \mathcal{F}$.

\begin{lemma}\label{lem3.1}
Let $M \in\mathcal{A}$ with $M \neq [0,\infty]$. Then there are a uniquely determined non-empty countable system $S=\{[b_{k}, d_{k}] \subseteq [0, \infty]\mid k\in K\}$ of closed intervals of a positive length which satisfy that for all $[b_{k}, d_{k}], [b_{l}, d_{l}]\in S, [b_{k}, d_{k}]\cap[b_{l}, d_{l}]=\emptyset $ or $ [b_{k}, d_{k}]\cap[b_{l}, d_{l}]=\{d_{k}\} $ when $d_{k}\leq b_{l}$, and a uniquely determined non-empty countable set $C=\{c_{k} \in [0, \infty] \mid k \in \overline{K}\}$ such that
$[b_{k}, d_{k}]\cap C\in\{\{b_{k}\}, \{d_{k}\}, \{b_{k}, d_{k}\}\}$ for all $k \in K$ and
\begin{equation*}
M= \{c_k\in [0,\infty]\mid k\in \overline{K}\}\cup \left([0,\infty]\setminus \left(\bigcup_{k\in K}[b_k, d_k] \right)\right)
\end{equation*}
where $|K|\leq|\overline{K}|$.
\end{lemma}
\begin{proof}The proof will be split into two parts.\\
$\mathbf{Part\, I}$. The existence of both $S$ and $C$ will be proved.

From $M \in\mathcal{A}$ with $M\neq [0,\infty]$, there is a function $f\in\mathcal{F}$ with Ran$(f)=M$. Take $S=\{[f(x^{-}),f(x)]\mid x\in [0, 1], f(x^{-})<f(x)\}\cup \{[f(x),f(x^{+})]\mid x \in [0, 1], f(x)< f(x^{+})\}$ and
$C=M_1\cup M_2\cup M_3$ where $M_1=\{f(x^{+})\mid x\in [0, 1], f(x^{+})\in M\mbox{ and }f(x^{-})< f(x^{+})\}$, $M_2=\{f(x^{-})\mid x \in [0, 1], f(x^{-})\in M \mbox{ and } f(x^{-})< f(x^{+})\}$ and $M_3=\{f(x)\mid x \in [0, 1], f(x^{-})< f(x^{+})\}\setminus (M_1\cup M_2)$.
 Clearly the set $S$ is countable. Let $K$ be a countable index set with $|K| =|S|$ and rewrite $S$ as $S=\{[b_{k}, d_{k}]\subseteq [0, \infty]\mid k \in K\}$ in which either $b_{k} = f(k^{-}), d_{k} = f(k)$ or $b_{k}=f(k), d_{k}=f(k^+)$ for each $k\in K$.
 The set $C$ is also countable.
Let $\overline{K}$ be a countable index set with $|\overline{K}|=|C|$. It is obvious that $|K|\leq |\overline{K}|$. Rearrange $C$ from small to large and assign its every element
an index $k\in \overline{K}$, for instance, the $k$-th element in the rearranged set $C$ is denoted by $c_{k}$ with $k\in \overline{K}$. This follows that $C= \{c_{k} \in [0, \infty] \mid k \in \overline{K}\}$. Note that for every $k\in K$, there is an $i \in \overline{K}$ such that $b_{k}=f(k)=c_{i}$ when $ f(k^-)\in M$ and $f(k)=f(k^-)$,
 there are two elements $s,t \in \overline{K}$ with $s + 1=t$ satisfying that $d_{k} = f(k)=c_{t}$ and $b_{k}=f(k^{-})= c_{s}$ when $f(k^-)\in M$ and $f(k)\neq f(k^-)$,
and there is a $j \in \overline{K}$ such that $d_{k}=f(k)= c_{j}$ when $ f(k^-)\notin M$.
Obviously, $S$ and $C$ have all required properties, respectively.\\
$\mathbf{Part\,II}$. The uniqueness of both $S$ and $C$ will be proved.

Suppose that a system $S_{1} = \{[u_{l}, v_{l}] \subseteq [0, \infty]\mid l \in K_{1} \}$ and a set $C_{1}= \{c_{l} \in [0, \infty] \mid l \in\overline{K}_{1}\}$ also have all required properties. Then we shall prove that $S=S_{1}$. We just prove $S\subseteq S_{1}$, and the case $S\supseteq S_{1}$ is completely analogous.

Let $[b_{k}, d_{k}]$ be an arbitrary element of $S$. Choosing $t\in [b_{k}, d_{k}]$
such that $t\notin M$, this means that there is a $[u_{l}, v_{l}]\in S_{1}$ such that $t\in[u_{l}, v_{l}]$ where
$[u_{l}, v_{l}] \cap M \in\{\{u_{l}\},\{v_{l}\}, \{u_{l}, v_{l}\}\}$.
 Next, we shall prove that $[b_{k}, d_{k}]=[u_{l}, v_{l}]$.

Supposing that $d_{k}\leq u_{l}< v_{l}$, this contradicts the fact that $t\in [b_{k}, d_{k}]$ and $t\in[u_{l}, v_{l}]$ with $t\notin M$.
Now, assume $u_{l}<d_{k}< v_{l}$. If $ d_{k}\in M$, then $d_{k}\in[u_{l}, v_{l}]\cap M$, contrary to $[u_{l}, v_{l}]\cap M\in\{\{u_{l}\},\{v_{l}\}, \{u_{l}, v_{l}\}\}$. If $d_{k}\notin M$, then $d_{k}$ is an accumulation point of
$M$ from the right. This follows that the set $[u_{l}, v_{l}]\cap M$ is infinite, a contradiction since $[u_{l}, v_{l}]\cap M\in\{\{u_{l}\},\{v_{l}\}, \{u_{l}, v_{l}\}\}$. Therefore, $u_{l}< v_{l}\leq d_{k}$.

In a completely analogous way, we have $b_{k}\leq u_{l}< v_{l}$.

Consequently, $b_{k}\leq u_{l}< v_{l}\leq d_{k}$.

Below, we shall prove that $v_{l}=d_{k}$, and the proof of $b_{k}=u_{l}$ is analogous. Suppose that $v_{l}<d_{k}$. Then $b_{k}< v_{l}< d_{k}$. If $ v_{l}\in M$, then $v_{l}\in[b_{k}, d_{k}]\cap M$, contrary to $[b_{k}, d_{k}]\cap M\in\{\{b_{k}\}, \{d_{k}\}, \{b_{k}, d_{k}\}\}$. If $v_{l}\notin M$, then $v_{l}$ is an accumulation point of
$M$ from the right. This follows that the set $[b_{k}, d_{k}]\cap M$ is infinite, a contradiction since $[b_{k}, d_{k}]\cap M\in\{\{b_{k}\}, \{d_{k}\}, \{b_{k}, d_{k}\}\}$. Therefore, $v_{l}=d_{k}$.


We finally come to $[b_{k}, d_{k}]=[u_{l}, v_{l}]$. From the arbitrariness of $[b_{k}, d_{k}]$, we have $S\subseteq S_{1}$.

Therefore, $S=S_{1}$. This follows that $C=C_{1} $. In particular, both $S$ and $C$ are independent of a choice of $f$.
\end{proof}

\begin{definition}\label{def3.1}
\emph{Let $M\in \mathcal{A}$. A pair $(S,C)$ is said to be associated with $M\neq [0,\infty]$ if $S=\{[b_{k}, d_{k}] \subseteq [0, \infty]\mid k\in K \}$ is a non-empty system of closed intervals of a positive length which satisfy that for all $[b_{k}, d_{k}],[b_{l}, d_{l}]\in S, [b_{k}, d_{k}]\cap[b_{l}, d_{l}]=\emptyset $ or $[b_{k}, d_{k}]\cap[b_{l}, d_{l}]=\{d_{k}\}$ when $d_{k}\leq b_{l}$, and $C=\{c_{k}\in [0, \infty]\mid k\in \overline{K}\}$ is a non-empty countable set such that
$[b_{k}, d_{k}]\cap C\in\{\{b_{k}\}, \{d_{k}\}, \{b_{k}, d_{k}\}\}$ for all $k\in K$, and
\begin{equation*}
M= \{c_k\in [0,\infty]\mid k\in \overline{K}\}\cup \left([0, \infty]\setminus \left(\bigcup_{k\in K}[b_k, d_k] \right)\right).
\end{equation*}
A pair $(S,C)$ is said to be associated with $M=[0,\infty]$ if $S=\{[\infty,\infty]\}$ and $C=\{\infty\}$.}
\end{definition}

We shall briefly write $(S,C)=(\{[b_k, d_k]\mid k\in K\}, \{c_k \mid k\in \overline{K}\})$ instead of $(S,C)=(\{[b_k, d_k] \subseteq [0, \infty]\mid k\in K\}, \{c_k\in [0,\infty]\mid k\in \overline{K}\})$.

The following remark is immediately from Definition \ref{def3.1}.
\begin{remark}\label{rem3.1} \emph{Let $M\in \mathcal{A}$, $(S,C)=(\{[b_k, d_k]\mid k\in K\}, \{c_k \mid k\in \overline{K}\})$ be associated with $M$ and $f\in \mathcal{F}$ with \mbox{Ran}$(f)=M$.
\renewcommand{\labelenumi}{(\roman{enumi})}
\begin{enumerate}
\item Let $x\in M\setminus C$. Then there is a $t\in (0,1)$ such that $f(t)=x$. It is easily seen that $f$ is continuous
at $t$ and $x$ is an accumulation point of $M$. In particular, $M\setminus C\subseteq \mbox{Acc}(M)$.
\item For all $a,b\in M\setminus C$, if $a<b$ then $(a,b)\cap(M\setminus C)\neq\emptyset$.
\item If $[b_k, d_k]\cap M=\{b_k\}$ for a certain $k\in K$, then $d_{k}$ is an accumulation point of $M$ from the right.
\item If $[b_k, d_k]\cap M=\{d_k\}$ for a certain $k\in K$, then $b_{k}$ is an accumulation point of $M$ from the left.
\item If $[b_k, d_k]\cap M=\{b_k,d_k\}$ for a certain $k\in K$, then there exist two elements $x,y\in[0,1]$ such that $f(x)=f(y)=b_k$.
\end{enumerate}}
\end{remark}

\begin{example}\label{exp3.1}
\renewcommand{\labelenumi}{(\roman{enumi})}
\emph{\begin{enumerate}
\item Let the function $f_1:[0,1]\rightarrow [0,\infty]$ be defined by\begin{equation*}
 f_1(x)=\begin{cases}
4x & \hbox{if }\ x\in[0,\frac{1}{2}),\\
2 & \hbox{if }\ x\in[\frac{1}{2},\frac{3}{4}),\\
9 & \hbox{if }\ x\in[\frac{3}{4},\frac{7}{8}),\\
10 & \hbox{if }\ x=\frac{7}{8},\\
\frac{2}{1-x} & \hbox{if }\ x\in(\frac{7}{8},1),\\
\infty & x=1
\end{cases}
\end{equation*}
Then the pair $(\{[2,9],[9,10],[10,16]\},\{9,10\})$ is associated with $[0,2)\cup \{9\}\cup (10)\cup [16,\infty]\in \mathcal{A}$.
\item Let the function $f_2:[0,1]\rightarrow [0,\infty]$ be defined by\begin{equation*}
 f_2(x)=\begin{cases}
8x & \hbox{if }\ x\in[0,\frac{1}{4}),\\
3 & \hbox{if }\ x\in[\frac{1}{4},\frac{1}{2}),\\
5 & \hbox{if }\ x\in[\frac{1}{2},\frac{3}{4}).\\
8x & \hbox{if }\ x\in[\frac{3}{4},\frac{7}{8}].\\
7 & \hbox{if }\ x\in(\frac{7}{8},1].
\end{cases}
\end{equation*}
Then the pair $(\{[2,3],[3,5], [5,6],[7,\infty]\},\{3,5,6,7\})$ is associated with $[0,2)\cup \{3\}\cup \{5\}\cup [6,7]\in \mathcal{A}$.
\item Let the function $f_3:[0,1]\rightarrow [0,\infty]$ be defined by \begin{equation*}
 f_3(x)=\begin{cases}
1 & \hbox{if }\ x\in[0,\frac{1}{10}),\\
10x & \hbox{if }\ x\in[\frac{1}{10},\frac{2}{5}),\\
8 & \hbox{if }\ x\in[\frac{2}{5},\frac{3}{5}),\\
10 & \hbox{if }\ x=\frac{3}{5},\\
20x& \hbox{if }\ x\in(\frac{3}{5},1].
\end{cases}
\end{equation*}
Then the pair $(\{[0,1],[4,8],[8,10],[10,12],[20,\infty],\{1,8,10,20\})$ is associated with $[1,4)\cup \{8\}\cup \{10\}\cup (12,20]\in \mathcal{A}$.
\end{enumerate}}
\end{example}

\section{An operation on $\mbox{Ran}(f)$ and its properties}

This section first defines an operation $\otimes$ on $\mbox{Ran}(f)$ with $f \in \mathcal{F}$, and then gives a necessary and sufficient condition for the operation $\otimes$ being associative that exactly depends on $\mbox{Ran}(f)$.

Let $x\in [0,\infty]$. Then by convention, $\sup\emptyset=0$ and $\inf\emptyset=\infty$.

\begin{definition}\label{def4.1}
\emph{Let $M\in\mathcal{A}$. Define a function $G_{M}:[0,\infty]\rightarrow M$ by
\begin{equation*}
 G_{M}(x)= \mbox{max}\{M\cap\{\sup([0,x]\cap M), \inf([x,\infty]\cap M)\}\}
\end{equation*}
for all $x\in[0,\infty]$.}
\end{definition}

The next proposition shows the relationship between $M$ and $G_{M}$.
\begin{proposition}\label{prop4.1}
Let $M\in \mathcal{A}$ and $(S,C)=(\{[b_k, d_k]\mid k\in K\}, \{c_k \mid k\in \overline{K}\})$ be associated with $M$. Then for all $x\in[0,\infty]$ and $k\in K$,
\renewcommand{\labelenumi}{(\roman{enumi})}
\begin{enumerate}
\item $G_{M}(x)=x$ if and only if $x\in M$.
\item if $x\notin M$ and $ x<f(0)$ then $G_{M}(x)=f(0)$.
\item if $x\notin M$ and $f(0)<x< f(1)$, then $G_{M}(x)=d_k$ if and only if $x\in [b_k, d_k]$ and $ d_k\in M$.
\item if $x\notin M$ and $f(0)<x< f(1)$, then $G_{M}(x)=b_k$ if and only if $x\in [b_k, d_k]$ and $ d_k\notin M$.
\item if $x\notin M$ and $f(1)<x$ then $G_{M}(x)=f(1)$.
\item $G_{M}$ is a non-decreasing function.
\end{enumerate}
\end{proposition}
\begin{proof}
Statements (i), (ii), (v) and (vi) are immediately from Definition \ref{def4.1}.

(iii) Let $x\notin M$ and $f(0)<x< f(1)$. If $G_{M}(x)=d_k$, then from the definition of $G_{M}$, $x\in [b_k, d_k]$ and $d_k\in M$. The converse implication is obvious.

(iv) Let $x\notin M$ and $f(0)<x< f(1)$. If $G_{M}(x)=b_k$, then from the definition of $G_{M}$, $x\in [b_k, d_k]$ and $d_k\notin M$. Conversely, if $x\in [b_k, d_k]$ and $d_k\notin M$, then $b_k\in M$, and we get $\sup([0,x]\cap M)=b_k$ and $\inf([x,\infty]\cap M)=d_k$, therefore, $G_{M}(x)=b_k$.
\end{proof}

\begin{example}\label{exp3.2}
\emph{In Example \ref{exp3.1}},
\emph{\renewcommand{\labelenumi}{(\roman{enumi})}
\begin{enumerate}
\item \begin{equation*}
  G_{M}(x)=\begin{cases}
9 & \rm{if }\ x\in[2,9],\\
10 & \rm{if }\ x\in(9,16],\\
x & \rm{otherwise}.
\end{cases}
\end{equation*}
\item \begin{equation*}
  G_{M}(x)=\begin{cases}
3 & \rm{if }\ x\in(2,3],\\
5 & \rm{if }\ x\in(3,5],\\
6 & \rm{if }\ x\in(5,6],\\
7 & \rm{if }\ x\in[7,\infty],\\
x & \rm{otherwise}.
\end{cases}
\end{equation*}
\item \begin{equation*}
 G_{M}(x)=\begin{cases}
1 & \rm{if }\ x\in[0,1],\\
8 & \rm{if }\ x\in[4,8],\\
10 & \rm{if }\ x\in(8,12],\\
20 & \rm{if }\ x\in[20,\infty],\\
x & \rm{otherwise}.
\end{cases}
\end{equation*}
\end{enumerate}}
\end{example}

In the forthcoming content, let $F:[0,\infty]^2\rightarrow [0,\infty]$ be an associative function.

\begin{definition}\label{def4.2}
\emph{Let $M\in\mathcal{A}$ and $G_{M}$ be determined by $M$. Define an operation $\otimes:M^2\rightarrow M $ by
\begin{equation*}
x\otimes y=G_{M}(F(x,y)).
\end{equation*}}
\end{definition}

\begin{example}\label{exp3.3}
\emph{In Example \ref{exp3.2}},
\emph{\renewcommand{\labelenumi}{(\roman{enumi})}
\begin{enumerate}
\item \begin{equation*}
 x\otimes y=\begin{cases}
9& \rm{if }\ F(x,y)\in[2,9],\\
10 & \rm{if }\ F(x,y)\in(9,16],\\
F(x,y) & \rm{otherwise}.
\end{cases}
\end{equation*}
\item \begin{equation*}
 x\otimes y=\begin{cases}
3 & \rm{if }\ F(x,y)\in(2,3],\\
5& \rm{if }\ F(x,y)\in(3,5],\\
6& \rm{if }\ F(x,y)\in(5,6],\\
7 & \rm{if }\ F(x,y)\in[7,\infty],\\
F(x,y) & \rm{otherwise}.
\end{cases}
\end{equation*}
\item \begin{equation*}
 x\otimes y=\begin{cases}
1 & \rm{if }\ F(x,y)\in[0,1],\\
8 & \rm{if }\ F(x,y)\in[4,8],\\
10 & \rm{if }\ F(x,y)\in(8,12],\\
20 & \rm{if }\ F(x,y)\in[20,\infty],\\
F(x,y) & \rm{otherwise}.
\end{cases}
\end{equation*}
\end{enumerate}}
\end{example}

\begin{proposition}\label{prop4.2}
Let $M\in\mathcal{A}$ and $(S, C)=(\{[b_k, d_k]\mid k\in K\}, \{c_k \mid k\in \overline{K}\})$ be associated with $M$. Then for all $x,y\in M$ and $k\in K$,
\renewcommand{\labelenumi}{(\roman{enumi})}
\begin{enumerate}
\item $x\otimes y=F(x,y)$ if and only if $F(x,y)\in M$.
\item if $F(x,y)\notin M$ and $F(x,y)<f(0)$ then $x\otimes y=f(0)$.
\item if $F(x,y)\notin M$ and $f(0)< F(x,y)< f(1)$, then $x\otimes y=d_k$ if and only if $F(x,y)\in [b_k, d_k]$ and $d_k\in M$.
\item if $F(x,y)\notin M$ and $f(0)< F(x,y)< f(1)$, then $x\otimes y=b_k$ if and only if $F(x,y)\in [b_k, d_k]$ and $d_k \notin M$
\item if $F(x,y)\notin M$ and $f(1)< F(x,y)$ then $x\otimes y=f(1)$.
\item $\otimes$ is a non-decreasing function.
\end{enumerate}
\end{proposition}
\begin{proof}
It is an immediate consequence of Proposition \ref{prop4.1} and Definition \ref{def4.2}.
\end{proof}

\begin{lemma}\label{lem4.1}
Let $M\in\mathcal{A}$ and $(S,C)=(\{[b_k, d_k]\mid k\in K\}, \{c_k \mid k\in \overline{K}\})$ be associated with $M$. Then $G_{M}(x)=f(f^{(-1)}(x))$ for all $x\in[0,\infty]$.
\end{lemma}
\begin{proof}
If $x\in \mbox{Ran}(f)$, then $f(f^{(-1)}(x))=x$ since $f\in \mathcal{F}$.

If $x\notin \mbox{Ran}(f)$ then there are three cases as follows.

(i) If $f(0)<x< f(1)$, then there is a $k\in K$ such that $x\in [b_k, d_k]$.
 Obviously, if $d_{k}\in M$ then
$f(f^{(-1)}(x))=d_{k}$;
if $d_{k}\notin M$ then
$f(f^{(-1)}(x))=b_{k}$.
Therefore, by Proposition \ref{prop4.1} (iii) and (iv), we get $G_{M}(x)=f(f^{(-1)}(x))$.

(ii) If $x< f(0)$, then $G_{M}(x)=f(0)$ and $f(f^{(-1)}(x))=f(\sup\{y\in [0,1]\mid f(y)<x\})=f(\sup \emptyset)=f(0)$.

(iii) If $f(1)< x$, then $G_{M}(x)=f(1)$ and $f(f^{(-1)}(x))=f(\sup\{y\in [0,1]\mid f(y)<x\})=f(1)$.

In summary, $G_{M}(x)=f(f^{(-1)}(x))$ for all $x\in[0,\infty]$.
\end{proof}

It is obvious that if a function $f$ is constant then $G_{M}(x)=f(f^{(-1)}(x))$ for all $x\in[0,\infty]$. Therefore, in the following we further suppose that all functions $f\in \mathcal{F}$ are non-constant. In particular, we have the following remark.
\begin{remark}\label{rem4.1} \emph{Let $f\in\mathcal{F}$. In Lemma \ref{lem4.1}, the condition that  $f(x)=f(x^{+})$ when $f(x^{+})\in \mbox{Ran}(f)$ cannot be deleted when $f$ isn't strictly monotone.
For instance, let the function $f:[0,1]\rightarrow [0,\infty]$ be defined by
\begin{equation*}
 f(x)=\begin{cases}
4x & \hbox{if }\ x\in[0,\frac{1}{4}],\\
2 & \hbox{if }\ x\in(\frac{1}{4},\frac{1}{2}],\\
8x & \hbox{if }\ x\in(\frac{1}{2},1].
\end{cases}
\end{equation*}
It is clear that $f(\frac{1}{4}^{+})\in \mbox{Ran}(f)$ but $f(\frac{1}{4})\neq f(\frac{1}{4}^{+})$ and $G_{M}(2)=2>1=f(f^{(-1)}(2))$, i.e., $G_{M}\neq f(f^{(-1)})$.}
\end{remark}

Let $M\in\mathcal{A}$ and $(S,C)=(\{[b_k, d_k]\mid k\in K\}, \{c_k \mid k\in \overline{K}\})$ be associated with $M$.
Denote
$$E=\{c\mid \mbox{ there are an } x_0\in[0,1]\mbox{ and }\varepsilon >0\mbox{ such that }f|_{[x_0,x_0+\varepsilon]}=c\},$$
$$N=\{\min \{x\in[0,1]\mid f(x)=y\}\mid y\in E\},~J=\{x\in [0,1]\mid f(x)\in M\setminus E \},$$
$$B=N\cup J.$$
One easily check that $f^{(-1)}(x) \in B$ for all $x\in [0,\infty]$.

We further give the following definition.
\begin{definition}\label{def4.3} \emph{Let $f\in\mathcal{F}$.
Define a function $f^*:B\rightarrow [0,\infty]$ by
\begin{equation*}
  f^*(x)= f(x)
\end{equation*}
for all $x\in B$ and a two-place function $F_0:B^2\rightarrow B$ by
 \begin{equation*}
F_0(x,y)=f^{(-1)}(F(f^*(x),f^*(y)))
\end{equation*}}
\end{definition}
for all $x,y\in B$, respectively.

Then from Definition \ref{def4.3}, one can check the following remark.
\begin{remark}\label{rem4.2} \emph{Let $f \in \mathcal{F}$. Then
\renewcommand{\labelenumi}{(\roman{enumi})}
\begin{enumerate}
\item $f^*$ is a strictly increasing function.
\item $f^{(-1)}(f^*(x))=x$ for all $x\in B$.
\end{enumerate}}
\end{remark}

Furthermore, we have the following proposition.
\begin{proposition}\label{prop4.3}
Let $f \in \mathcal{F}$. Then
 \begin{equation*}
x\otimes y=f^*(F_0(f^{(-1)}(x),f^{(-1)}(y)))
\end{equation*}
for all $x,y\in M$ and
 \begin{equation*}
F_0(x,y)=f^{(-1)}(f^*(x)\otimes f^*(y))
\end{equation*}
for all $x,y\in B$. Moreover, $F_0$ is associative if and only if $\otimes$ is associative.
\end{proposition}
\begin{proof} For any $x,y\in M$, it follows from Definitions \ref{def4.2} and \ref{def4.3}, Lemma \ref{lem4.1} and Remark \ref{rem4.2} that
\begin{eqnarray*}
x\otimes y&=&G_{M}(F(x,y))\\
&=&f(f^{(-1)}(F(x,y)))\\
&=&f^*(f^{(-1)}(F(x,y)))\\
&=&f^*(f^{(-1)}(F(f^*(f^{(-1)}(x)),f^*(f^{(-1)}(y)))))\\
&=&f^*(F_0(f^{(-1)}(x),f^{(-1)}(y))).
\end{eqnarray*}
There are two elements $u,v\in B$ such that $f^*(u)=x$, $f^*(v)=y$ since $x,y\in M$. Then from Remark \ref{rem4.2},
\begin{equation*}
f^*(u)\otimes f^*(v)=x\otimes y=f^*(F_0(f^{(-1)}(x),f^{(-1)}(y)))=f^*(F_0(u,v)),
\end{equation*}
which deduces $f^{(-1)}(f^*(u)\otimes f^*(v))=f^{(-1)}(f^*(F_0(u,v))=F_0(u,v)$. Clearly, this follows that $F_0$ is associative if and only if $\otimes$ is associative.
\end{proof}

\begin{lemma}\label{lem4.2}
Let $f \in \mathcal{F}$ and $T:[0,1]^2\rightarrow [0,1]$ be a function defined by Eq.(\ref{eq:5}). Then, for any $x,y \in [0,1]$, there exist two elements $m,n\in B$ such that $f^*(m)=f(x)$, $f^*(n)=f(y)$ and $T(x,y)=F_0(m,n)$. In particular, $T(x,y)=F_0(x,y)$ for all $x,y\in B$.
\end{lemma}
\begin{proof}
If $x,y\in B$, then from Definition \ref{def4.3}, we are sure that $f^*(x)=f(x)$, $f^*(y)=f(y)$. If $x\notin B$, then $f(x)\in E$. Let $m=\min \{s\in[0,1]\mid f(s)=f(x), f(x)\in E\}$. Obviously, $m\in B$ and $f^*(m)=f(x)$. Analogously, if $y\notin B$ then there exists an $n\in B$ such that $f^*(n)=f(y)$. Therefore, from Definition \ref{def4.3}, we get
\begin{eqnarray*}
T(x,y)&=&f^{(-1)}(F(f(x),f(y)))\\
&=&f^{(-1)}(F(f^*(m),f^*(n)))\\
&=&F_0(m,n).
\end{eqnarray*}
\end{proof}

\begin{proposition}\label{prop4.4}Let $f \in \mathcal{F}$ and $T:[0,1]^2\rightarrow [0,1]$ be a function defined by Eq.(\ref{eq:5}).
Then $F_0$ is associative if and only if $T$ is associative.
\end{proposition}
\begin{proof}
Suppose that $F_0$ is associative. Let us prove that $T(T(x,y),z)=T(x,T(y,z))$ for all $x,y,z\in[0,1]$.
 For any $x,y,z\in [0,1]$, by Lemma \ref{lem4.2} there exist three elements $m,n,t\in B$ such that $f^*(m)=f(x)$, $f^*(n)=f(y)$, $f^*(t)=f(z)$, respectively, and $T(x,y)=F_0(m,n)$, $T(x,z)=F_0(m,t)$, $T(y,z)=F_0(n,t)$, respectively. Then
 \begin{eqnarray*}
T(T(x,y),z)&=&T(F_0(m,n),z)\\
&=&F_0(F_0(m,n),t)\\
&=&F_0(m,F_0(n,t))\\
&=&F_0(m,T(y,z))\\
&=&T(x,T(y,z)).
\end{eqnarray*}

Conversely, if $T$ is associative, i,e., $T(T(x,y),z)=T(x,T(y,z)$ for all $x,y,z \in [0,1]$, then by Eq.(\ref{eq:5}), $$f^{(-1)}(F(f\circ f^{(-1)} (F(f(x),f(y))),f(z)))=f^{(-1)}(F(f(x),f\circ f^{(-1)} (F(f(x),f(y))))).$$
 Therefore, for all $x,y,z\in B$, $$f^{(-1)}(F(f^*\circ f^{(-1)} (F(f^*(x),f^*(y))),f^*(z)))=f^{(-1)}(F(f^*(x),f^*\circ f^{(-1)} (F(f^*(x),f^*(y))))).$$
 Finally, from Definition \ref{def4.3}, $F_0(F_0(x,y),z)=F_0(x,F_0(y,z)$, i.e., $F_0$ is associative.
\end{proof}

The following is an immediate consequence of Propositions \ref{prop4.3} and \ref{prop4.4}, which means that the associativity of $T$ defined by Eq.(\ref{eq:5}) relies only on $\mbox{Ran}(f)$.
\begin{theorem}\label{theorem4.1}
Let $f \in \mathcal{F}$ and $T:[0,1]^2\rightarrow [0,1]$ be a function defined by Eq.(\ref{eq:5}).
Then $T$ is associative if and only if $\otimes$ is associative.
\end{theorem}

\section{Associativity of the function $T$ defined by Eq.(\ref{eq:5})}

This section shows two necessary and sufficient conditions for the function $T$ defined by Eq.(\ref{eq:5}) being associative.

In the rest of this article, for every $k\in K$, if there are $x,y \in M$ such that $F(x,y)\in[b_{k},d_{k}]\setminus \{c_{k}\}$ where $\{c_{k}\}=M\cap[b_{k},d_{k}]$ then we use the symbol $x\otimes y=a_{k}$. From Proposition \ref{prop4.2}, we have $a_{k}=b_{k}$ or $a_k=d_{k}$.

Let $M\subseteq [0,\infty]$. Define $O(M) = \bigcup_{x,y\in M}[\min\{x,y\}, \max\{x,y\}]$ when $M\neq \emptyset$, and $O(M)=\emptyset$ when $M=\emptyset$. It is clear that $O(M)=\emptyset$ if and only if $M=\emptyset$. Denote $F(A,c)=\{F(x,c)\mid x\in A\}$, $F(c,A)=\{F(c,x)\mid x\in A\}$ and $F(\emptyset,A)=F(A,\emptyset)=\emptyset$.
For every $k\in K$, let $I_{k}=O(\{a_{k}\}\cup\{z\in[b_{k},d_{k}]\setminus \{c_{k}\}\mid \mbox{ there are two elemenets }x,y \in M \mbox{ such that } F(x,y)=z\})$.

In what follows we shall explore another characterization that the function $T$ defined by Eq.(\ref{eq:5}) is associative. To achieve the goal, the following lemma plays a central role.
\begin{lemma}\label{lem5.2}
Let $M\in \mathcal{A}$ and $(S,C)=(\{[b_k, d_k]\mid k\in K\}, \{c_k \mid k\in \overline{K}\})$ be associated with $M$. For any $x,y\in[0,\infty]$, if $[\min\{x,y\},\max\{x,y\}]\cap M= \emptyset$ then $G_M(x)= G_M(y)$.
\end{lemma}
\begin{proof}Supposing $G_M(x)\neq G_M(y)$, from Proposition \ref{prop4.1} we get $G_M(\min\{x,y\})=s<t=G_M(\max\{x,y\})$ with $s,t\in M$. Obviously, $[s,t]\cap M\neq\emptyset$. Choose $w\in [s,t]\cap M$. Then $w=G_M(w)$. This follows that $G_M(\min\{x,y\})\leq G_M(w)\leq G_M(\max\{x,y\})$ and $\min\{x,y\}\leq w\leq \max\{x,y\}$ since $G_M$ is non-decreasing. Hence, $[\min\{x,y\},\max\{x,y\}]\cap M\neq\emptyset$, a contradiction.
\end{proof}

The following proposition gives a sufficient condition that the operation $\otimes$ on $M$ is associative.
\begin{proposition}\label{Propo5.1}
Let $M\in \mathcal{A}$ and $(S,C)=(\{[b_k, d_k]\mid k\in K\}, \{c_k \mid k\in \overline{K}\})$ be associated with $M$. If $F(\cup_{k\in K}I_{k}, M)\cap M= \emptyset$ and
$F(M, \cup_{k\in K}I_{k})\cap M= \emptyset$, then the operation $\otimes$ on $M$ is associative.
\end{proposition}
\begin{proof} Supposing that $F(\cup_{k\in K}I_{k}, M)\cap M=\emptyset$, we shall prove that for any $x_{1},x_{2},x_{3}\in M$, $G_{M}(F(F(x_{1},x_{2}),x_{3}))=(x_{1}\otimes x_{2})\otimes x_{3}$. Indeed, from Definition \ref{def4.2} we have $x_{1}\otimes x_{2}=G_{M}(F(x_{1},x_{2}))$.

We distinguish two cases as follows.

(i) If $F(x_{1},x_{2})\in M$, then $x_{1}\otimes x_{2}=F(x_{1},x_{2})$, which implies
\begin{align}
 G_{M}(F(F(x_{1},x_{2}),x_{3})) &= G_{M}(F(x_{1}\otimes x_{2},x_{3})) \nonumber \\
 &= (x_{1}\otimes x_{2})\otimes x_{3}.\nonumber
\end{align}

(ii) If $F(x_{1},x_{2})\notin M$, then there is a $k \in K$ such that $F(x_{1},x_{2})\in [b_{k}, d_{k}]\setminus{\{c_{k}\}}$, where $\{c_{k}\} = M\cap[b_{k}, d_{k}]$. Thus from $F(\cup_{k\in K}I_{k}, M)\cap M= \emptyset$,
  $$F(I_{k},x_{3})\cap M= \emptyset.$$
This means
$$F([\min\{a_{k}, F(x_{1},x_{2})\}, \max\{a_{k}, F(x_{1},x_{2})\}], x_{3})\cap M= \emptyset$$
 since $[\min\{a_{k}, F(x_{1},x_{2})\}, \max\{a_{k}, F(x_{1},x_{2})\}] \in I_{k}$.
Further, from Lemma \ref{lem5.2},
$$G_{M}(F(F(x_{1},x_{2}),x_{3}))= G_{M}(F(a_{k},x_{3})),$$
which deduces
$$G_{M}(F(a_{k},x_{3}))=G_{M}(F(x_{1}\otimes x_{2},x_{3}))= (x_{1}\otimes x_{2})\otimes x_{3}$$
since $x_{1}\otimes x_{2}=a_{k}$. Hence $G_{M}(F(F(x_{1},x_{2})),x_{3})=(x_{1}\otimes x_{2})\otimes x_{3}$.

$G_{M}(F(x_{1},F(x_{2},x_{3}))=x_{1}\otimes (x_{2}\otimes x_{3})$ is completely analogous when $(F(M, \cup_{k\in K}I_{k}))\cap M= \emptyset$.

Therefore, $(x_{1}\otimes x_{2})\otimes x_{3}=x_{1}\otimes (x_{2}\otimes x_{3})$ since $F$ is associative.
\end{proof}

Proposition \ref{Propo5.1} implies the following corollary.
\begin{corollary}\label{cor5.1}
Let $M\in \mathcal{A}$ and $(S,C)=(\{[b_k, d_k]\mid k\in K\}, \{c_k \mid k\in \overline{K}\})$ be associated with $M$. If $F$ is commutative and $F(\cup_{k\in K}I_{k}, M)\cap M= \emptyset$,
then the operation $\otimes$ on $M$ is associative.
\end{corollary}

In particular, we have the following corollary.
\begin{corollary}\label{cor5.11}
Let $M\in \mathcal{A}$ and $(S,C)=(\{[b_k, d_k]\mid k\in K\}, \{c_k \mid k\in \overline{K}\})$ be associated with $M$. If $F(M, M)\subseteq M$,
then the operation $\otimes$ on $M$ is associative.
\end{corollary}
\begin{proof}
In fact, if $F(M, M) \subseteq M$ then $\cup_{k\in K}I_{k}= \emptyset$, implying $F(\cup_{k\in K}I_{k}, M)= \emptyset$ and $(F(M, \cup_{k\in K}I_{k}))= \emptyset$.
Further, by Proposition \ref{Propo5.1}, the operation $\otimes$ on $M$ is associative.
\end{proof}

Proposition \ref{Propo5.1} implies the following remark.
\begin{remark}\label{rem5.1}\emph{Let $f \in \mathcal{F}$, $T:[0,1]^2\rightarrow [0,1]$ be a function defined by Eq.(\ref{eq:5}) in which $F:[0,\infty]^2\rightarrow [0,\infty]$ is a non-decreasing and commutative function.
\begin{enumerate}\renewcommand{\labelenumi}{(\roman{enumi})}
\item Let $0$ be a neutral element of $F$, $(F(\cup_{k\in K}I_{k}, M))\cap M= \emptyset$ and $F(M, \cup_{k\in K}I_{k})\cap M= \emptyset$. Then $T$ is associative. Meanwhile, $T$ isn't necessary a t-superconorm since $T(x,0)=f^{(-1)}(F(f(x),f(0)))\geq f^{(-1)}(f(x))$ and $f^{(-1)}(f(x))\leq x$ for all $x\in [0,1]$. However, if $f$ is further a strictly increasing function with $f(0)=0$, then $T$ is a t-conorm.
\item  Let $\infty$ be a neutral element of $F$, $(F(\cup_{k\in K}I_{k}, M))\cap M= \emptyset$ and $F(M, \cup_{k\in K}I_{k})\cap M= \emptyset$. Then $T$ is a t-subnorm. Moreover, $T$ isn't necessary a t-norm since $T(x,1)=f^{(-1)}(F(f(x),f(1)))\leq f^{(-1)}(f(x))\leq x$ for all $x\in [0,1]$. However, if $f$ is further a strictly increasing function with $f(1)=\infty$, then $T$ is a t-norm. Another way is to slightly modify the function $T$ as
 \begin{equation}\label{eq:1.3}
T(x,y)=\left\{
 \begin{array}{ll}
  f^{(-1)}(F(f(x),f(y))) & \hbox{if }(x,y)\in [0,1)^{2}, \\
  \min\{x,y\}& \hbox{otherwise}
 \end{array}
\right.
\end{equation}
 for all $x,y\in [0,1]$. Then one may check that $T$ is a t-norm.
\end{enumerate}}
\end{remark}

Generally, the converse of Proposition \ref{Propo5.1} isn't true.
\begin{example}\label{exap5.1}
\emph{Let $F(x,y)=x+y$ and the function $f:[0,1]\rightarrow [0,\infty]$ be defined by
\begin{equation*}
 f(x)=x \ \ \hbox{for all }\ x\in[0,1].
\end{equation*}
Then from Eq.(\ref{eq:5}),
\begin{equation*}
T(x,y)=\left\{
 \begin{array}{ll}
  x+y & \hbox{if } x+y<1, \\
  1 & \hbox{otherwise.}
 \end{array}
\right.
\end{equation*}
It is easy to check that $T$ is associative, $M=[0,1]$, $\cup_{k\in K}I_{k}=[1,2]$ and $F(\cup_{k\in K}I_{k}, M)=[1,3]$. Thus $F(\cup_{k\in K}I_{k}, M) \cap M\neq\emptyset$.}
\end{example}

A binary function $F:[0, \infty]^2\rightarrow [0, \infty]$ satisfies the cancellation law if  $F(x,y)=F(x,z)$ implies $x=0$ ($x=\infty$) or $y=z$.
Fortunately, we have the following proposition.
\begin{proposition}\label{prop5.5}
 Let $M\in \mathcal{A}$ and $(S,C)=(\{[b_k, d_k]\mid k\in K\}, \{c_k \mid k\in \overline{K}\})$ be associated with $M$. If $F $ is cancellative and satisfies $F(M,M\setminus C)\subseteq M\setminus C$, then the following are equivalent:
\renewcommand{\labelenumi}{(\roman{enumi})}
\begin{enumerate}
\item The operation $\otimes$ is associative.
\item $F(\cup_{k\in K}I_{k}, M)\cap M=\emptyset$ and $F(M, \cup_{k\in K}I_{k})\cap M= \emptyset$.
\end{enumerate}
\end{proposition}
\begin{proof}
From Proposition \ref{Propo5.1}, it is enough to prove that (ii) implies (i). We shall prove $F(M, M) \subseteq M$. Indeed, if $F(M, M) \nsubseteq M$, then there exist two elements $x, y \in M$ such that $F(x , y)\notin M $. Thus
$F(x , y) \in [b_{k}, d_{k}]\setminus{\{c_{k}\}}$ for some $k \in K$ where $\{c_{k}\}=[b_{k}, d_{k}] \cap M$. Now, let $z \in M\setminus C$.
Then by $F((M,M\setminus C)) \subseteq M\setminus C$, we have $F(y,z) \in M\setminus C$. This means $x \otimes (y \otimes z) = x \otimes F(y , z) = G_{M}(F(x, F(y,z))$. On the other hand $(x \otimes y) \otimes z = a_{k}\otimes z = G_{M}(F(a_{k}, z))$. Thereby, the associativity of $\otimes$ on $M$ yields $G_{M}(F(x, F(y,z)))=G_{M}(F(a_{k}, z))$. This follows $G_{M}(F(x, F(y,z)))= F(a_{k}, z)$ since $F(a_{k}, z)\in M\setminus C$. Therefore, by Proposition \ref{prop4.1} (i), $F(x, F(y,z))= F(a_{k}, z)$ since $F(x, F(y,z))\in F(M,M\setminus C)\subseteq M\setminus C$, which results in $F(F(x,y),z))= F(a_{k}, z)$ since $F $ associative. So that $F(x , y)=a_{k}$ since $F $ is cancellative, contrary to $F(x , y) \in [b_{k}, d_{k}]\setminus{\{c_{k}\}}$. Therefore, $F(M, M) \subseteq M$, then $\cup_{k\in K}I_{k}= \emptyset$, implying $(F(\cup_{k\in K}I_{k}, M))\cap M= \emptyset$ and $F(M, \cup_{k\in K}I_{k})\cap M= \emptyset$.
\end{proof}

Then Theorem \ref{theorem4.1} and Proposition \ref{prop5.5} imply the following characterization theorem that the function $T$ defined by Eq.(\ref{eq:5}) is associative.
\begin{theorem}\label{theorem5.10}
 Let $M\in \mathcal{A}$ and $(S,C)=(\{[b_k, d_k]\mid k\in K\}, \{c_k \mid k\in \overline{K}\})$ be associated with $M$. If $F $ is cancellative and satisfies $F(M,M\setminus C)\subseteq M\setminus C$, then the following are equivalent:
\renewcommand{\labelenumi}{(\roman{enumi})}
\begin{enumerate}
\item $T$ is associative.
\item $F(\cup_{k\in K}I_{k}, M)\cap M=\emptyset$ and $F(M, \cup_{k\in K}I_{k})\cap M= \emptyset$.
\end{enumerate}
\end{theorem}

Notice that in Theorem \ref{theorem5.10} the condition that $F $ is cancellative and satisfies $F(M,M\setminus C)\subseteq M\setminus C$ cannot be deleted generally.
\begin{example}\label{exap5.21}\emph{Let $M\in \mathcal{A}$ and $(S,C)=(\{[b_k, d_k]\mid k\in K\}, \{c_k \mid k\in \overline{K}\})$ be associated with $M$.}\\
\emph{(1) Let $F(x,y)=x+y-xy$ and the function $f:[0,1]\rightarrow [0,\infty]$ be defined by
\begin{equation*}
 f(x)=\left\{
 \begin{array}{ll}
  0 & \hbox{if } x\in[0,\frac{1}{2}), \\
  x & \hbox{if } x\in[\frac{1}{2},1].
 \end{array}
\right.
\end{equation*}
Then from Eq.(\ref{eq:5}),
\begin{equation*}
T(x,y)=\left\{
 \begin{array}{ll}
 0 & \hbox{if } (x,y)\in[0,\frac{1}{2})\times[0,\frac{1}{2}),\\
 x+y-xy & \hbox{if } (x,y)\in[\frac{1}{2},1]\times[\frac{1}{2},1],\\
 \max\{x,y\} & \hbox{otherwise.}
 \end{array}
\right.
\end{equation*}
 One may check that $F$ isn't cancellative, $M=\{0\}\cup[\frac{1}{2},1]$, $F(M,M\setminus C)=(\frac{1}{2},1)\subseteq M\setminus C$ and $F(\cup_{k\in K}I_{k}, M)=\{0\}$. Thus, $F(\cup_{k\in K}I_{k}, M) \cap M\neq\emptyset$. However, $T$ is associative.}\\
 \emph{(2) Let $F(x,y)=x+y$ and the function $f:[0,1]\rightarrow [0,\infty]$ be defined by
\begin{equation*}
 f(x)=\left\{
 \begin{array}{ll}
  x & \hbox{if } x\in[0,\frac{1}{2}), \\
  1 & \hbox{if } x\in[\frac{1}{2},1].
 \end{array}
\right.
\end{equation*}
Then from Eq.(\ref{eq:5}),
\begin{equation*}
T(x,y)=\left\{
 \begin{array}{ll}
 \min\{x,y\} & \hbox{if } x+y<1,\\
  1 & \hbox{otherwise.}
 \end{array}
\right.
\end{equation*}
 One may check that $F$ is cancellative, $M=[0,\frac{1}{2})\cup\{1\}$, $F(M,M\setminus C)=[0,\frac{3}{2})\nsubseteq M\setminus C$ and $(F(\cup_{k\in K}I_{k}, M))=[\frac{1}{2},3)$. Thus, $F(\cup_{k\in K}I_{k}, M) \cap M\neq\emptyset$. However, $T$ is associative.}
\end{example}

Below, we shall develop another characterization that the function $T$ defined by Eq.(\ref{eq:5}) is associative. We first need the following definition.
\begin{definition}\label{def5.1}
\emph{Let $M\in \mathcal{A}$ and $(S,C)=(\{[b_k, d_k]\mid k\in K\}, \{c_k \mid k\in \overline{K}\})$ be associated with $M$. Define for all $y\in M$, and $k,l\in K$, $M_{k}^{y}=\{x\in M \mid F(x,y)\in [b_k, d_k]\setminus\{c_k\}\},~ ~M_{y}^{k}=\{x\in M \mid F(y,x)\in [b_k, d_k]\setminus\{c_k\} \},$ $M^{y}=\{x\in M \mid F(x,y)\in M\setminus C\},~M_{y}=\{x\in M \mid F(y,x)\in M\setminus C\}$,
$$H_{k}^{y}=\{(x_1,x_2)\in M_{k}^{y}\times M_{y} \mid F(F(x_1,y),x_2)\neq F(a_k,x_2), \hbox{if } M_{k}^{y}\neq\emptyset,  M_{y}\neq\emptyset\},$$
$$H_{y}^{k}=\{(x_1,x_2)\in M^{y} \times  M_{y}^{k}\mid F(x_1,F(y,x_2))\neq F(x_1, a_k), \hbox{if } M_{y}^{k}\neq\emptyset,  M^{y}\neq\emptyset\},$$
 $$H_{k,l}^{y}=\{(x_1,x_2)\in M_{k}^{y}\times M_{l}^{y} \mid F(a_k,x_2)\neq F(x_1,a_l), \hbox{if } M_{k}^{y}\neq\emptyset,  M_{l}^{y}\neq\emptyset\},$$
 \begin{equation*}
J_{k,l}^{y}=\begin{cases}
O(F(M_{k}^{y},a_l)\cup F(a_k,M_{y}^{l})), & \hbox{if }\ M_{k}^{y}\neq\emptyset, M_{y}^{l}\neq\emptyset,\\
\emptyset, & \hbox{otherwise,}
\end{cases}
\end{equation*}
  $I_{k}^{y}=O(\{a_k\}\cup F(M_{k}^{y},y))$ and $I_{y}^{k}=O(\{a_k\}\cup F(y,M_{y}^{k}))$.
Put $\mathfrak{J_{1}}(M)=\bigcup_{y\in M}\bigcup_{k\in K}F(I_{k}^{y},M_{y}),~ \mathfrak{J_{2}}(M)=\bigcup_{y\in M}\bigcup_{k\in K}F(M^{y},I_{y}^{k})$, $\mathfrak{J_{3}}(M)=\bigcup_{y\in M}\bigcup_{k,l\in K}J_{k,l}^{y}\mbox{ and}$,
$$\mathfrak{J}(M)=\mathfrak{J_{1}}(M)\cup\mathfrak{J_{2}}(M)\cup\mathfrak{J_{3}}(M).$$}
\end{definition}

Then we have the following sufficient condition that the function $T$ defined by Eq.(\ref{eq:5}) is associative.
\begin{proposition}\label{prop5.10}
Let $M\in \mathcal{A}$ and $(S,C)=(\{[b_k, d_k]\mid k\in K\}, \{c_k \mid k\in \overline{K}\})$ be associated with $M$, $\otimes$ be the operation on $M$. If $\mathfrak{J}(M)\cap M=\emptyset$
then the operation $\otimes$ on $M$ is associative.
\end{proposition}
\begin{proof}
 Suppose that the operation $\otimes$ isn't associative. Then there exist three elements $x,y,z\in M$ such that $(x\otimes y)\otimes z\neq x\otimes(y\otimes z)$. We claim that $F(x,y)\notin M\setminus C$ or $F(y,z)\notin M\setminus C$. Otherwise, from Definition \ref{def4.1}, $F(x,y)\in M\setminus C$ and $F(y,z)\in M \setminus C$ would imply $(x\otimes y)\otimes z=G_M(F(F(x,y),z))= G_M(F(x,F(y,z)))=x\otimes(y\otimes z)$, a contradiction. The following proof is completed by three steps.

(i) Let $F(x,y)\notin M\setminus C$ and $F(y,z)\in M\setminus C$. Then $y\otimes z=F(y,z)$ and there exists a $k\in K$ such that $F(x,y)\in [b_k,d_k]\setminus\{c_k\}$. Thus $x\otimes y=a_k$. It follows from Definition \ref{def4.1} that $G_M(F(a_k,z))=(x\otimes y)\otimes z\neq x\otimes(y\otimes z)=G_M(F(x,F(y,z)))$. On the other hand, by the associativity of $F$, we have $G_M(F(x,F(y,z)))=G_M(F(F(x,y),z))$. Thus $G_M(F(a_k,z))\neq G_M(F(F(x,y),z))$. Therefore, by Lemma \ref{lem5.2}, $F([\min\{a_k,F(x,y)\},\max\{a_k,F(x,y)\}],z)\cap M\neq\emptyset $. Obviously, $z\in M_{y}$, and $x\in M_{k}^{y}$. Thus by $I_{k}^{y}=O(\{a_k\}\cup F(M_{k}^{y},y))$, $F([\min\{a_k,F(x,y)\},\max\{a_k,F(x,y)\}],z)\subseteq F(I_{k}^{y},M_{y})$, this infers $F(I_{k}^{y},M_{y})\cap M \neq\emptyset$.

(ii) Let $F(x,y)\in M\setminus C$ and $F(y,z)\notin M\setminus C$. In a completely analogous way, we have $F(M^{y},I_{y}^{k})\cap M \neq\emptyset$.

(iii) Let $F(x,y)\notin M\setminus C$ and $F(y,z)\notin M\setminus C$. Then from $F(x,y)\notin M\setminus C$, there exists a $k\in K$ such that $F(x,y)\in [b_k,d_k]\setminus\{c_k\}$. Thus $x\otimes y=a_k$. Similarly, there exists an $l\in K$ such that $F(y,z)\in [b_l,d_l]\setminus\{c_l\}$ since $F(y,z)\notin M\setminus C$. Hence $y\otimes z=a_l$. Therefore, $G_M(F(a_k,z))=(x\otimes y)\otimes z\neq x\otimes(y\otimes z)=G_M(F(x,a_l))$. It follows from Lemma \ref{lem5.2} that $[\min\{F(a_k,z),F(x,a_l)\}, \max\{F(a_k,z),F(x,a_l)\}]\cap M \neq \emptyset$, which together with $x\in M_{k}^{y}$, $F(x,a_l)\in F(M_{k}^{y},a_l)$, $z\in M_{y}^{l}$ and $F(a_k,z)\in F(a_k,M_{y}^{l})$ yields that $[\min\{F(a_k,z),F(x,a_l)\},\max\{F(a_k,z),F(x,a_l)\}]\subseteq J_{k,l}^{y}$. Therefore, $J_{k,l}^{y}\cap M \neq \emptyset$.

From (i), (ii) and (iii), we deduce that $\mathfrak{J}(M)\cap M\neq\emptyset$.
\end{proof}

Generally, the converse of Proposition \ref{prop5.10} isn't true. For instance, in Example \ref{exap5.1}, $M=[0,1]$, $\mathfrak{J_{1}}(M)=\mathfrak{J_{2}}(M)=[1,3]$ and $\mathfrak{J_{3}}(M)=[1,2]$, thus, $\mathfrak{J}(M) \cap M\neq\emptyset$.

To obtain the necessary condition that the function $T$ defined by Eq.(\ref{eq:5}) is associative, we further suppose that $F:[0,\infty]^2\rightarrow [0,\infty]$ is a monotone and associative function with neutral element in $[0,1]$. Then we have the following two lemmas.
\begin{lemma}\label{lem5.3}
Let $M\in \mathcal{A}$ and $(S,C)=(\{[b_k, d_k]\mid k\in K\}, \{c_k \mid k\in \overline{K}\})$ be associated with $M$. Let $M_1, M_2\subseteq [0,\infty]$ be two non-empty sets and $c\in[0,\infty]$. If there exist $u\in M_1$ and $v\in M_2$ such that $F(u,c)\neq F(v,c)$ and $F(c,u)\neq F(c,v)$,
then\\
(1) $F(O(M_1\cup M_2),c)\cap (M\setminus C)=\emptyset$ if and only if there exist $x_{1}\in M_1$ and $x_{2}\in M_2$ such that
$$[\min\{F(x_{1},c),F(x_{2},c)\}, \max\{F(x_{1},c),F(x_{2},c)\}]\cap (M\setminus C)=\emptyset.$$
(2) $F(c,O(M_1\cup M_2))\cap (M\setminus C)=\emptyset$ if and only if there exist $x_{1}\in M_1$ and $x_{2}\in M_2$ such that
 $$[\min\{F(c,x_{1}),F(c,x_{2})\}, \max\{F(c,x_{1}),F(c,x_{2})\}]\cap (M\setminus C)=\emptyset.$$
 \end{lemma}
\begin{proof}
(1) We only provide the proof of the statement when $F$ is non-decreasing, and the statement is analogous when $F$ is non-increasing.

If $F(O(M_1\cup M_2),c)\cap (M\setminus C)=\emptyset$, then it is easy to see that there exist an $x_{1}\in M_1$ and an $x_{2}\in M_2$ such that
$$[\min\{F(x_{1},c),F(x_{2},c)\}, \max\{F(x_{1},c),F(x_{2},c)\}]\cap (M\setminus C)=\emptyset.$$

Conversely, suppose that $F(O(M_1\cup M_2),c)\cap (M\setminus C)\neq\emptyset$. Now, let $a\in F(O(M_1\cup M_2),c)\cap (M\setminus\mathcal{C})$. Because there are $u\in M_1$ and $v\in M_2$ such that $F(u,c)\neq F(v,c)$, there exist two elements $m,n\in M_1\cup M_2$ with $m<n$ such that $F(m,c)<F(n,c)$ and $a\in [F(m,c),F(n,c)]$ since $F$ is non-decreasing. It is easy to see that the point $a$ is an accumulation point of $M$ from the left and right. Hence,  there exist two elements $m,n\in M_1\cup M_2$ with $m<n$ such that
$F([m, n], c)\cap (M\setminus C)\neq\emptyset.$
The following proof is split into three cases.
\renewcommand{\labelenumi}{(\roman{enumi})}
\begin{enumerate}
\item In the case that exactly one of $\{m,n\}$ is contained in $M_1$ and the other is contained in $M_2$, it is obvious that
$$[\min\{F(x_{1},c),F(x_{2},c)\}, \max\{F(x_{1},c),F(x_{2},c)\}]\cap (M\setminus C)\neq\emptyset.$$
\item In the case that $m,n\in M_1$ with $m<n$. Choose $t\in M_2$.
\renewcommand{\labelenumi}{(\roman{enumi})}
\begin{enumerate}
\item If $t\leq m$ then put $x_{1}=t$ and $x_{2} =n$.
\item If $t\geq n$ then put $x_{1}=m$ and $x_{2} =t$.
\item If $m<t<n$, then $F([m, t], c)\cap (M\setminus C)\neq\emptyset$ or $F([t, n], c)\cap (M\setminus C)\neq\emptyset$. If $F([m, t], c)\cap (M\setminus C)\neq\emptyset$, then put $x_{1}=m$ and $x_{2} =t$. If $F([t, n], c)\cap (M\setminus C)\neq\emptyset$, then put $x_{2}=n$ and $x_{1} =t$.
\end{enumerate}
In any case of (a),(b) and (c), we get $$[\min\{F(x_{1},c),F(x_{2},c)\}, \max\{F(x_{1},c),F(x_{2},c)\}]\cap (M\setminus C)\neq\emptyset.$$
\item The case $m,n\in M_2$ is completely analogous to (ii).
\end{enumerate}
(2) The proof is analogous to the proof of the statement (1).
\end{proof}

\begin{lemma}\label{lem5.1}
Let $M\in \mathcal{A}$ and $(S,C)=(\{[b_k, d_k]\mid k\in K\}, \{c_k \mid k\in \overline{K}\})$ be associated with $M$. If $G_M(x)= G_M(y)$ for any $x\neq y\in[0,\infty]$ then $[\min\{x,y\},\max\{x,y\}]\cap (M\setminus C)= \emptyset$.
\end{lemma}
\begin{proof} Suppose that $[\min\{x,y\},\max\{x,y\}]\cap (M\setminus C)\neq\emptyset$. Then say $x<y$, thus we have
$[x,y]\cap (M\setminus C)\neq\emptyset$. Let {$a\in [x,y]\cap (M\setminus C)$}. Then $x\leq a\leq y$ and $G_M(a)=a$, which deduce that $G_M(x)\leq G_M(a)=a\leq G_M(y)$ since $G_M$ is a non-decreasing function.
If $y\in M\setminus C$, then from Proposition \ref{prop4.1}, $G_M(y)=y$, it follows that $G_M(x)< y= G_M(y)$, a contradiction.
If $y\notin M\setminus C$, then there exists a $k\in K$ such that $y\in[b_k,d_k]$ and $G_M(y)=b_k$ or $G_M(y)=d_{k}$. So that $a<b_k$ since $a\in M\setminus C$ and $a<y$, hence $G_M(a)<G_M(y)$. Therefore, $G_M(x)\leq G_M(a)< G_M(y)$, a contradiction.
\end{proof}

Note that under the condition of Lemma \ref{lem5.1} we cannot get $[\min\{x,y\},\max\{x,y\}]\cap M= \emptyset$ generally.
\begin{example}\label{exap5.20}
\emph{ Let  $f:[0,1]\rightarrow [0,\infty]$ be defined by
\begin{equation*}
 f(x)=\left\{
 \begin{array}{ll}
  x & \hbox{if } x\in[0,\frac{1}{4}), \\
  \frac{1}{2} & \hbox{if } x=\frac{1}{2}, \\
  2x & \hbox{if } x\in(\frac{1}{2},1].
 \end{array}
\right.
\end{equation*}
It is easy to check that $M=[0,\frac{1}{4})\cup\{\frac{1}{2}\}\cup(1,2]$ and $G_M(\frac{3}{8})=\frac{1}{2}= G_M(\frac{3}{4})$. Thus $[\min\{x,y\},\max\{x,y\}]=[\frac{3}{8}, \frac{3}{4}]$ and $[\frac{3}{8},\frac{3}{4}]\cap (M\setminus C)= \emptyset$. However, $[\frac{3}{8},\frac{3}{4}]\cap M \neq \emptyset$. Therefore, the converse of Lemma \ref{lem5.2} isn't true generally.}
\end{example}

Let $M\in \mathcal{A}$, $(\mathcal{S},\mathcal{C})=(\{[b_k, d_k]\mid k\in K\}, \{c_k \mid k\in K\})$ be associated with $M$. For all $k,l\in K, y\in M$, write

$(C_1)$ either $H_{k}^{y} = \emptyset$ or $F(I_{k}^{y}, M_{y}) \cap (M \setminus C) = \emptyset$,

$(C_2)$ either $H_{y}^{k} = \emptyset$ or $F(M^{y},I_{y}^{k}) \cap (M \setminus C) = \emptyset$,

$(C_3)$ either $H_{k,l}^{y} = \emptyset$ or $J_{k,l}^{y} \cap (M \setminus C) = \emptyset$,

Conditions $(C_1)$, $(C_2)$  and $(C_3)$ are called an $F$-condition of $M$.

We further have the following proposition that is almost the converse of Proposition \ref{prop5.10}.
\begin{proposition}\label{prop5.11}
Let $M\in \mathcal{A}$ and $(S,C)=(\{[b_k, d_k]\mid k\in K\}, \{c_k \mid k\in \overline{K}\})$ be associated with $M$.
If the operation $\otimes$ on $M$ is associative then the $F$-condition of $M$ holds.
\end{proposition}
\begin{proof}Suppose the $F$-condition of $M$ does not hold. Then there exist a $y\in M$ and two elements $k,l\in K$ such that $H_{k}^{y} \neq \emptyset$ and $F(I_{k}^{y},M_{y})\cap (M\setminus C) \neq\emptyset$,  $H_{y}^{k} = \emptyset$ and $F(M^{y},I_{y}^{k})\cap (M\setminus C) \neq\emptyset$ or, $H_{k,l}^{y} \neq \emptyset$ and $J_{k,l}^{y}\cap (M\setminus C)\neq \emptyset$. We distinguish three cases as follows.

(i) If $H_{k}^{y} \neq \emptyset$ and $F(I_{k}^{y},M_{y})\cap (M\setminus C) \neq\emptyset$, then there exists a $z\in M_{y}$ such that $F(I_{k}^{y},z)\cap (M\setminus C) \neq\emptyset$. Thus by the definition of $I_{k}^{y}$, $F(M_{k}^{y},y)\neq\emptyset$. Because of $H_{k}^{y}\neq\emptyset$, Applying Lemma \ref{lem5.3}, there are a $u\in C$ and a $v\in F(M_{k}^{y},y)$ such that
$[\min\{F(u,z),F(v,z)\},\max\{F(u,z),F(v,z)\}]\cap (M\setminus C)\neq\emptyset.$
Furthermore, there is an $x\in M_{k}^{y}$ such that $F(x,y)=v$ since $v\in F(M_{k}^{y},y)$. Thereby, there are a $u\in C$ and an $x\in M_{k}^{y}$ such that
$[\min\{F(a_k,z),F(F(x,y),z)\},\max\{F(a_k,z),F(F(x,y),z)\}]\cap (M\setminus C)\neq\emptyset.$
Consequently, from Lemma \ref{lem5.1}, we get $G_M(F(a_k,z))\neq G_M(F(F(x,y),z))$.
On the other hand, from $x\in M_{k}^{y}$ we have $F(x,y)\in [b_k,d_k]\setminus\{c_k\}$. Thus $x\otimes y=a_k$. From $z\in M_{y}$, we get $F(y,z)\in M\setminus C$. This follows $y\otimes z=F(y,z)$.
Therefore, $(x\otimes y)\otimes z=G_M(F(a_k,z))\neq G_M(F(F(x,y),z))=G_M(F(x,F(y,z)))= x\otimes(y\otimes z)$.

(ii) If $H_{y}^{k} = \emptyset$ and $F(M^{y},I_{y}^{k})\cap (M\setminus C) \neq\emptyset$, then in complete analogy to (i), one can check that $(x\otimes y)\otimes z\neq x\otimes(y\otimes z)$.

(iii) If $J_{k,l}^{y}\cap (M\setminus C)\neq \emptyset$, then $J_{k,l}^{y}\neq \emptyset$. By the definition of $J_{k,l}^{y}$, $F(O(F(M_{k}^{y},a_l)\cup F(a_k,M_{y}^{l})),e)\cap (M\setminus C)\neq \emptyset$ where $e$ is a neutral element of $F$, $F(M_{k}^{y},a_l)\neq\emptyset$ and $F(a_k,M_{y}^{l})\neq\emptyset$. Because of $H_{k}^{y}\neq\emptyset$, Using Lemma \ref{lem5.3}, there are a $u\in F(a_k,M_{y}^{l})$ and a $v\in F(M_{k}^{y},a_l)$ such that $$[\min\{F(u,e),F(v,e)\},\max\{F(u,e),F(v,e)\}]\cap (M\setminus C)\neq\emptyset.$$
Since $u\in F(a_k,M_{y}^{l})$ and $v\in F(M_{k}^{y},a_l)$, there are an $x\in M_{k}^{y}$ and a $z\in M_{y}^{l}$ such that $u=F(a_k,z)$, $v=F(x,a_l)$. Therefore, there are an $x\in M_{k}^{y}$ and a $z\in M_{y}^{l}$ such that
$$[\min\{F(a_k,z),F(x,a_l)\},\max\{F(a_k,z),F(x,a_l)\}]\cap (M\setminus C)\neq\emptyset$$
since $e$ is a neutral element of $F$. Furthermore, by Lemma \ref{lem5.1} we have $G_M(F(a_k,z))\neq G_M(F(x,a_l))$.
On the other hand, from $x\in M_{k}^{y}$ we have $F(x,y)\in [b_k,d_k]\setminus\{c_k\}$. Thus $x\otimes y=a_k$.
From $z\in M_{y}^{l}$, we get $F(y,z)\in [b_l,d_l]\setminus\{c_l\}$. Thus $y\otimes z=a_l$.
We finally come to $(x\otimes y)\otimes z=G_M(F(a_k,z))\neq G_M(F(x,a_l))= x\otimes(y\otimes z)$.
\end{proof}

Generally, the converse of Proposition \ref{prop5.11} isn't true.

\begin{example}\label{exap5.2}
\emph{(i) Let $F(x,y)=x+y$ and the function $f:[0,1]\rightarrow [0,\infty]$ be defined by
\begin{equation*}
 f(x)=\left\{
 \begin{array}{ll}
  x & \hbox{if } x\in[0,\frac{1}{4}), \\
  \frac{1}{4} & \hbox{if } x\in[\frac{1}{4},\frac{1}{2}), \\
  1 & \hbox{if } x\in[\frac{1}{2},1].
 \end{array}
\right.
\end{equation*}
Then from Eq.(\ref{eq:5}),
\begin{equation*}
T(x,y)=\left\{
 \begin{array}{ll}
  x+y & \hbox{if } x+y\leq\frac{1}{4}, \\
  \frac{1}{4} & \hbox{if } (x,y)\in\{0\}\times(\frac{1}{4},\frac{1}{2})\cup(\frac{1}{4},\frac{1}{2})\times\{0\},\\
  1 & \hbox{otherwise.}
 \end{array}
\right.
\end{equation*}
 One may check that $M=[0,\frac{1}{4})\cup\{1\}$, $M\setminus C=[0,\frac{1}{4})$ and $\mathfrak{J}(M)=[\frac{1}{4},\frac{3}{2}]$, thus, $\mathfrak{J}(M) \cap (M\setminus C)=\emptyset$. However, $T$ isn't associative.}\\
 \emph{(ii) Let $F(x,y)=x+y$ and the function $f:[0,1]\rightarrow [0,\infty]$ be defined by
\begin{equation*}
 f(x)=\left\{
 \begin{array}{ll}
  x & \hbox{if } x\in[0,1), \\
  2 & \hbox{if } x=1.
 \end{array}
\right.
\end{equation*}
Then from Eq.(\ref{eq:5}),
\begin{equation*}
T(x,y)=\left\{
 \begin{array}{ll}
  x+y & \hbox{if } x+y<1, \\
  1 & \hbox{otherwise.}
 \end{array}
\right.
\end{equation*}
 One may check that $M=[0,1)\cup\{2\}$, $M\setminus C=[0,1)$ and $\mathfrak{J}(M)=[1,5)$. Thus $\mathfrak{J}(M) \cap (M\setminus C)=\emptyset$ and $T$ is associative. However, $\mathfrak{J}(M) \cap M\neq \emptyset$.}
\end{example}

With Propositions \ref{prop5.10} and \ref{prop5.11}, we have the following one.
\begin{proposition}\label{prop5.12}
Let $M\in \mathcal{A}$ and $(S,C)=(\{[b_k, d_k]\mid k\in K\}, \{c_k \mid k\in \overline{K}\})$ be associated with $M$. Then $F(C, M)\cup F(M,C)\subseteq M\setminus C$ implies that the following are equivalent:\\
(1) The operation $\otimes$ is associative.\\
(2) the $F$-condition of $M$ holds.
\end{proposition}
\begin{proof}
(1) implying (2) is immediately from Proposition \ref{prop5.11}. Now, suppose that the operation $\otimes$ isn't associative. Then by the proof of Proposition \ref{prop5.10}, there exist three elements $x,y,z\in M$ such that either $F(x,y)\notin M\setminus C$ or $F(y,z)\notin M\setminus C$.

If $F(x,y)\notin M\setminus C$ and $F(y,z)\in M\setminus C$ then there is a $k\in K$ such that $F(x,y)\in [b_k,d_k]\setminus\{c_k\}$. Thus $x\otimes y=a_k$. Hence $G_M(F(a_k,z))\neq G_M(F(x,G_M(F(y,z))$, this together with $F(C, M)\subseteq M\setminus C$ deduces that $[\min\{F(a_k,z),F(x,F(y,z))\},\max\{F(a_k,z),F(x,F(y,z))\}]\cap (M\setminus C)\neq\emptyset $.
 Then, by the proof of Proposition \ref{prop5.10},  $F(a_k,z)\neq F(F(x,y),z)$ and $F(I_{k}^{y},M_{y})\cap (M\setminus C) \neq\emptyset$. Therefore, $H_{k}^{y} \neq \emptyset$ and $\mathfrak{J}(M)\cap (M\setminus C)\neq\emptyset$.

If $F(x,y)\in M\setminus C$ and $F(y,z)\notin M\setminus C$ then, in complete analogous, one can check that $H_{y}^{k} \neq \emptyset$ and $F(M^{y},I_{y}^{k})\cap (M\setminus C) \neq\emptyset$.

If $F(x,y)\notin M\setminus C$ and $F(y,z)\notin M\setminus C$, then there exists an $l\in K$ such that $F(y,z)\in [b_l,d_l]\setminus\{c_l\}$. Hence $y\otimes z=a_l$. Thus by the proof of Proposition \ref{prop5.10}, $$ F(a_k,z)\neq F(x,a_l)\mbox{ and }[\min\{F(a_k,z),F(x,a_l)\},\max\{F(a_k,z),F(x,a_l)\}]\subseteq J_{k,l}^{y}.$$ Therefore, $H_{k,l}^{y} \neq \emptyset$ and $J_{k,l}^{y}\cap (M\setminus C) \neq \emptyset$. This follows that the $F$-condition of $M$ does not hold.
\end{proof}

Theorem \ref{theorem4.1} and Proposition \ref{prop5.12} imply the following characterization theorem that the function $T$ defined by Eq.(\ref{eq:5}) is associative.
\begin{theorem}\label{theorem5.12}
Let $M\in \mathcal{A}$ and $T:[0,1]^2\rightarrow [0,1]$ be a function defined by Eq.(\ref{eq:5}). If $F(C, M)\cup F(M,C)\subseteq M\setminus C$, then the function $T$ is associative if and only if the $F$-condition of $M$ holds.
\end{theorem}

In particular, if $F$ is cancellative, then for all $k,l\in K, y\in M$, $H_{k}^{y} =H_{y}^{k}=H_{k,l}^{y}\neq\emptyset$. Hence, we have the following corollary.
\begin{corollary}\label{co02.3.001}
	Let $M\in \mathcal{A}$ and $T:[0,1]^2\rightarrow [0,1]$ be a function defined by Eq.(\ref{eq:5}). If $F$ is cancellative and $F(C, M)\cup F(M,C)\subseteq M\setminus C$, then the function $T$ is associative if and only if $\mathfrak{J}(M)\cap (M\setminus C)=\emptyset$.
\end{corollary}

In general, the condition $F(C, M)\cup F(M,C)\subseteq M\setminus C$ in Theorem \ref{theorem5.12} cannot be dropped .
\begin{example}\label{exap5.3}
\emph{(1) Let $F(x,y)=\frac{xy}{2}$ for all $x,y \in [0,\infty]$ and the function $f:[0,1]\rightarrow [0,\infty]$ be defined by \begin{equation*}
 f(x)=\begin{cases}
x & \hbox{if }\ x\in[0, \frac{1}{2}),\\
 \frac{1}{2} & \hbox{if }\ x\in[\frac{1}{2},1].
\end{cases}
\end{equation*}
From Eq.(\ref{eq:5}),
\begin{equation*}
T(x,y)=\frac{xy}{2} \hbox{ for all }\ x, y\in[0, 1].
\end{equation*}
One may check that $M=[0,\frac{1}{2}]$, $\mathfrak{J}(M)=\emptyset$ and $F(C, M)=F(M,C)=[0,\frac{1}{8}]\subseteq M\setminus C$. Moreover, $\mathfrak{J}(M) \cap (M\setminus C)=\emptyset$ and $T$ is an associative function.}\\
\emph{(2) Let $F(x,y)=x+y+1$ for all $x,y \in [0,\infty]$ and the function $f:[0,1]\rightarrow [0,\infty]$ be defined by \begin{equation*}
 f(x)=\begin{cases}
2x & \hbox{if }\ x\in[0, \frac{1}{2}),\\
 \frac{1}{1-x} & \hbox{if }\ x\in[\frac{1}{2},1),\\
 \infty & \hbox{if }\ x=1.
\end{cases}
\end{equation*}
Then from Eq.(\ref{eq:5}),
\begin{equation*}
T(x,y)=\left\{
 \begin{array}{ll}
  \frac{1}{2} & \hbox{if } x+y<\frac{1}{2}, \\
  \frac{2-x-y}{3-2x-2y+xy}& \hbox{otherwise.}
 \end{array}
\right.
\end{equation*}
One may check that $M=[0,1)\cup[2,\infty]$, $\mathfrak{J}(M)=[2,4)\cup (4, \infty]$ and $F(C, M)=F(M,C)=[3,\infty]\subseteq M\setminus C$. Thus, $\mathfrak{J}(M) \cap (M\setminus C)\neq\emptyset$ and $T$ isn't an associative function.}\\
\emph{(3) Let $F(x,y)=xy$ for all $x,y \in [0,\infty]$ and the function $f:[0,1]\rightarrow [0,\infty]$ be defined by \begin{equation*}
 f(x)=\begin{cases}
x & \hbox{if }\ x\in[0, \frac{1}{4}),\\
 1 & \hbox{if }\ x\in[\frac{1}{4},1].
\end{cases}
\end{equation*}
Then from Eq.(\ref{eq:5}),
\begin{equation*}
T(x,y)=\left\{
 \begin{array}{ll}
  xy & \hbox{if }(x,y)\in[0,\frac{1}{4})^{2}, \\
  1 & \hbox{if }(x,y)\in[\frac{1}{4},1]^{2}, \\
  \min\{x,y\} & \hbox{otherwise.}
 \end{array}
\right.
\end{equation*}
One may check that $M=[0,\frac{1}{4})\cup\{1\}$, $F(C, M)=F(M,C)=[0,\frac{1}{4})\cup\{1\}\nsubseteq M\setminus C$ and $\mathfrak{J}(M) \cap (M\setminus C)=[0,\frac{1}{4})$. However, $T$ is an associative function.}\\
\emph{(4) Let $F(x,y)=\max\{x,y\}$ for all $x,y \in [0,\infty]$ and the function $f:[0,1]\rightarrow [0,\infty]$ be defined by \begin{equation*}
 f(x)=\begin{cases}
x & \hbox{if }\ x\in[0, 1),\\
 2 & \hbox{if }\ x=1.
\end{cases}
\end{equation*}
Then from Eq.(\ref{eq:5}),
\begin{equation*}
T(x,y)=\left\{
 \begin{array}{ll}
  \max\{x,y\} & \hbox{if }(x,y)\in[0,1)^{2}, \\
  1& \hbox{otherwise.}
 \end{array}
\right.
\end{equation*}
One may check that $M=[0,1)\cup\{2\}$ and $F(C, M)=F(M,C)=\{2\}\nsubseteq M\setminus C$. However, $\mathfrak{J}(M) \cap (M\setminus C)=\emptyset$ and $T$ is an associative function.}
\end{example}

\section{Conclusions}

In this article, we proved that the function $T:[0,1]^2\rightarrow[0,1]$ defined by Eq.(\ref{eq:5}) is associative if and only if $F(\cup_{k\in K}I_{k}, M)\cap M=\emptyset$ and $F(M, \cup_{k\in K}I_{k})\cap M= \emptyset$ when $F $ is cancellative and satisfies $F(M,M\setminus C)\subseteq M\setminus C$ if and only if the $F$-condition of $M$ holds when $F(C, M)\cup F(M,C)\subseteq M\setminus C$, where $M=\mbox{Ran}(f)$ with $f: [0,1]\rightarrow [0,\infty]$ a non-decreasing function that satisfies either $f(x)=f(x^{+})$ when $f(x^{+})\in \mbox{Ran}(f)$ or $f(x)\neq f(y)$ for any $y\neq x$ when $f(x^{+})\notin \mbox{Ran}(f)$ for all $x\in[0,1]$. Example \ref{exap5.21} shows that, in general, one cannot delete the condition that $F $ is cancellative and satisfies $F(M,M\setminus C)\subseteq M\setminus C$. Meanwhile, Example \ref{exap5.3} reveals that the condition $F(C, M)\cup F(M,C)\subseteq M\setminus C$ cannot be moved generally. It is worth pointing out that one can verify that all the results obtained hold when $f$ is a non-increasing function. Therefore, our results are applicable for the monotone functions $f$ that satisfy either $f(x)=f(x^{+})$ when $f(x^{+})\in \mbox{Ran}(f)$ or $f(x)\neq f(y)$ for any $y\neq x$ when $f(x^{+})\notin \mbox{Ran}(f)$ for all $x\in[0,1]$. On the other hand, it is evident that our results are suitable for all strictly monotone functions and all monotone right continuous ones, respectively. Therefore, they generalize the work of Vicen\'{\i}k \cite{PV2005} and Zhang and Wang \cite{YM2024}. Comparing with Yao Ouyang\cite{YO2008}, our function $f$ needn't be strictly decreasing and satisfies the formula \eqref{eq:3}, which means that our results deeply generalize the work of Yao Ouyang\cite{YO2008}. An interesting problem is whether our results can be applied to all monotone functions or not.

\section*{Acknowledgments}
The authors thank the referees for their valuable comments and suggestions.


\begin{thebibliography}{99}
\bibitem{Abel}
N. H. Abel, Untersuchung der Function zweier unabhangig veranderlichen Grossen $x$ und $y$ wie $f(x,y)$, welche die Eigenschaft haben, dass $f(z,f(x,y))$ eine symmetrische Function von $x$, $y$ und $z$ ist, J. Reine Angew. Math. 1 (1826) 11-15.

\bibitem{CA2006}
 C. Alsina, M. J. Frank, B. Schweizer, Associative Functions: Triangular Norms and Copulas, World Scientific, Singapore, 2006.

\bibitem{SJ1999}
S. Jenei, Fibred triangular norms, Fuzzy Sets Syst. 103 (1999) 67-82.

\bibitem{EP2000}
 E. P. Klement, R. Mesiar, E. Pap, Triangular Norms, Kluwer Academic Publishers, Dordrecht, 2000.

\bibitem{CH1965}
C. H. Ling, Representation of associative functions, Publ. Math. Debrecen 12 (1965) 189-212.

\bibitem{AM2004}
A. Mesiarov\'{a}-Zem\'{a}nkov\'{a}, Continuous triangular subnorms, Fuzzy Sets Syst. 142 (2004) 75-83.

\bibitem{YO2007}
Y. Ouyang, On the construction of boundary weak triangular norms through additive generators, Nonlinear Anal. 66 (2007) 125-130.

\bibitem{YO2008}
Y. Ouyang, J. Y. Fang, Z. J. Zhao, A generalization of additive generator of triangular norms, Int. J. Approx. Reason. 49 (2008) 417-421.

\bibitem{BS1961}
B. Schweizer, A. Sklar, Associative functions and satistical triangle inequalities, Publ. Math. Debrecen 8 (1961) 169-186.
\bibitem{BS1963}
B. Schweizer, A. Sklar, Associative functions and abstract semigroups, Publ. Math. Debrecen 10 (1963) 69-81.
\bibitem{PV1998}
P. Vicen\'{\i}k, A note on generators of t-norms, BUSEFAL 75 (1998) 33-38.
\bibitem{PV1998b}
P. Vicen\'{\i}k, Additive generators and discontinuity, BUSEFAL 76 (1998) 25-28.

\bibitem{PV2005}
P. Vicen\'{\i}k, Additive generators of associative functions, Fuzzy Sets Syst. 153 (2005) 137-160.

\bibitem{PV2008}
P. Vicen\'{\i}k, Additive generators of border-continuous triangular norms, Fuzzy Sets Syst. 159 (2008) 1631-1645.

\bibitem{DZ2005}
D. X. Zhang, Triangular norms on partially ordered sets, Fuzzy Sets Syst. 153 (2005) 195-209.

\bibitem{YM2024}
Y. M. Zhang, X. P. Wang, Characterizations of monotone right continuous functions which generate associative functions, Fuzzy Sets Syst. 477 (2024) 108799.
\end{thebibliography}
\end{document}